A Series Expansion for Integral Powers of Arctangent

Michael Milgram[1], Consulting Physicist

Box 1484, Deep River, Ont. Canada. K0J 1P0

**Abstract**: A new series expansion for the function $[\arctan(x)/x]^n$ is developed, and some properties of the expansion coefficients are obtained.

**Introduction**: Although infinite series representations for the classical functions have been known for well over a century, I have found no indication from computerized literature searches that the expansion presented here (**Theorem 3**) is known. Certainly it does not appear in any of the standard reference works [**1**, **4**, **5**, **9** or **11**], nor can it be found in representative textbooks of yesteryear [**2**, especially Ch.8, **10**, particularly page 124].

However, this main result has considerable significance and application. As is well known, the Fourier transform of $\exp(-r)/r^2$ involves the arctangent function, and an expansion in powers of arctangent corresponds to the Neumann solution of integral equations for a range of problems of physical interest [**3**]. Thus, with reference to **Theorem 3**, it now becomes possible to convert hitherto intractable integrals involving higher powers of arctangent, into (slowly) converging power series. It also turns out that such series can be evaluated in closed form, thereby uncovering new, higher order solutions for a range of problems in physics. This is discussed elsewhere [**7**, **8**].

**Theorem 1**

$$\sum_{n=0}^{k}\frac{1}{k-n+\frac{1}{2}}\sum_{m=0}^{n}\frac{1}{m+\frac{1}{2}} = 2\sum_{n=0}^{k}\frac{1}{n+1}\sum_{m=0}^{n}\frac{1}{m+\frac{1}{2}} \qquad (1)$$

---

[1] mike@geometrics-unlimited.com

**Proof:**

Proceed by induction denoting the left- and right-hand sides by L (k) and R (k) respectively. **( 1 )** is obviously true for k=1 and has been verified for k=2. Then L(k)=R(k) for at least one value of "k", and it is sufficient to show that L(k+1)=R(k+1).

$$\begin{aligned}
L(k+1) &= \sum_{n=0}^{k+1} \frac{1}{k-n+\frac{3}{2}} \sum_{m=0}^{n} \frac{1}{m+\frac{1}{2}} \\
&= \sum_{n=0}^{k} \frac{1}{k-n+\frac{1}{2}} \sum_{m=0}^{n+1} \frac{1}{m+\frac{1}{2}} + \frac{2}{k+\frac{3}{2}} \\
&= \sum_{n=0}^{k} \frac{1}{k-n+\frac{1}{2}} \sum_{m=0}^{n} \frac{1}{m+\frac{1}{2}} + \frac{2}{k+\frac{3}{2}} + \sum_{n=0}^{k} (\frac{1}{k-n+\frac{1}{2}})(\frac{1}{n+\frac{3}{2}}) \\
&= R(k) + \frac{2}{k+\frac{3}{2}} + \sum_{n=0}^{k} (\frac{1}{k-n+\frac{1}{2}})(\frac{1}{n+\frac{3}{2}}) \\
&= R(k+1) - \frac{2}{(k+2)} \sum_{m=0}^{k+1} \frac{1}{m+\frac{1}{2}} + \frac{2}{k+\frac{3}{2}} + \sum_{n=0}^{k} (\frac{1}{k-n+\frac{1}{2}})(\frac{1}{n+\frac{3}{2}}) \\
&= R(k+1)
\end{aligned}$$

( 2 )

In the second last line, split the final term into partial fractions, reverse one of the resulting sums and simplify to obtain the final line. ▲

**Corollary**

**( 1 )** can be rewritten as a transformation between sums of digamma functions:

$$\sum_{l=1}^{k} \frac{\psi(l+\frac{1}{2})}{k-l+\frac{1}{2}} = 2\sum_{l=1}^{k} \frac{\psi(l+\frac{1}{2})}{l} + \psi(\tfrac{1}{2})[\psi(k+\tfrac{1}{2}) - \psi(\tfrac{1}{2}) + 2\psi(1) - 2\psi(k+1)] \qquad (3)$$

**Theorem 2**

$$\boxed{\begin{aligned}
&\sum_{m_n=0}^{k} \frac{1}{(k-m_n+\frac{1}{2})} \sum_{m_{n-1}=0}^{m_n} \frac{1}{(m_{n-1}+\frac{n-1}{2})} \sum_{m_{n-2}=0}^{m_{n-1}} \frac{1}{(m_{n-2}+\frac{n-2}{2})} \cdots \sum_{m_2=0}^{m_3} \frac{1}{(m_2+1)} \sum_{m_1=0}^{m_2} \frac{1}{(m_1+\frac{1}{2})} \\
&= n\sum_{m_n=0}^{k} \frac{1}{(m_n+\frac{n}{2})} \sum_{m_{n-1}=0}^{m_n} \frac{1}{(m_{n-1}+\frac{n-1}{2})} \sum_{m_{n-2}=0}^{m_{n-1}} \frac{1}{(m_{n-2}+\frac{n-2}{2})} \cdots \sum_{m_2=0}^{m_3} \frac{1}{(m_2+1)} \sum_{m_1=0}^{m_2} \frac{1}{(m_1+\frac{1}{2})}
\end{aligned}}$$

( 4 )

**Proof:**

Define coefficients

$$t_k(n-1) = \sum_{m_{n-1}=0}^{k} \frac{1}{(m_{n-1}+\frac{n-1}{2})} \sum_{m_{n-2}=0}^{m_{n-1}} \frac{1}{(m_{n-2}+\frac{n-2}{2})} \cdots \sum_{m_2=0}^{m_3} \frac{1}{(m_2+1)} \sum_{m_1=0}^{m_2} \frac{1}{(m_1+\frac{1}{2})} \quad (5)$$

or, recursively

$$t_k(n) = \sum_{m=0}^{k} \frac{1}{(m+\frac{n}{2})} t_m(n-1) \quad (6)$$

where $t_k(0) = 1$ and $t_k(J) = 0 \ \forall \ J < 0$. The $t_k(n)$ possess the following properties:

$$t_{k+1}(n) = t_k(n) + \frac{1}{(k+1+\frac{n}{2})} t_{k+1}(n-1) \quad (7)$$

$$t_0(n) = \frac{2}{n} t_0(n-1) \quad \text{or} \quad t_0(n) = \frac{2^n}{n!} \quad (8)$$

and are unexpectedly reminiscent of similar coefficients that independently arise in the theory of generalized exponential integrals of continuous order [6, Eq. 2.23]. In this notation, ( 4 ) reads

$$\sum_{m=0}^{k} \frac{1}{(k-m+\frac{1}{2})} t_m(n-1) = n \sum_{m=0}^{k} \frac{1}{(m+\frac{n}{2})} t_m(n-1) \quad (9)$$

Proceed by induction, noting that ( 9 ) is true for n=2 by virtue of **Theorem 1**. Therefore assume that the following (( 10 )) is true for some value of "n-1" and show that it is true for "n".

$$\sum_{m=0}^{k} \frac{1}{(k-m+\frac{1}{2})} t_m(n-2) = (n-1) \sum_{m=0}^{k} \frac{1}{(m+\frac{n-1}{2})} t_m(n-2) \quad (10)$$

Note that ( 9 ) is trivially true when k=0 for all values of "n". So, let the left- and right-hand sides of ( 9 ) be denoted by $L_k(n)$ and $R_k(n)$ respectively, and similarly assume $L_k(n)=R_k(n)$ for some value of "$k$". It is sufficient to show, for any "$n$", that $L_{k+1}(n)=R_{k+1}(n)$. Trivially, **( 9 )** reads

$$R_{k+1}(n) = R_k(n) + \frac{n}{(k+1+\frac{n}{2})} t_{k+1}(n-1) \qquad (11)$$

and

$$\begin{aligned}
L_{k+1}(n) &= \sum_{m=0}^{k+1} \frac{1}{k-m+\frac{3}{2}} t_m(n-1) = \sum_{m=-1}^{k} \frac{1}{k-m+\frac{1}{2}} t_{m+1}(n-1) \\
&= \sum_{m=0}^{k} \frac{1}{k-m+\frac{1}{2}} t_{m+1}(n-1) + \frac{t_0(n-1)}{(k+\frac{3}{2})} \\
&= \sum_{m=0}^{k} \frac{1}{k-m+\frac{1}{2}} t_m(n-1) + \sum_{m=0}^{k} \frac{1}{(k-m+\frac{1}{2})(m+1+\frac{n-1}{2})} t_{m+1}(n-2) \\
&\quad + \frac{t_0(n-1)}{(k+\frac{3}{2})}
\end{aligned} \qquad (12)$$

using **( 7 )**. The second term may be split into partial fractions, and the first term identified as $L_k(n)$, giving

$$\begin{aligned}
L_{k+1}(n) &= L_k(n) + \frac{t_0(n-1)}{(k+\frac{3}{2})} \\
&\quad + \frac{1}{(k+\frac{n+2}{2})} \{\sum_{m=0}^{k} \frac{1}{(m+1+\frac{n-1}{2})} t_{m+1}(n-2) + \sum_{m=0}^{k} \frac{1}{(k-m+\frac{1}{2})} t_{m+1}(n-2)\} \\
&= L_k(n) + \frac{t_0(n-1)}{(k+\frac{3}{2})} \\
&\quad + \frac{1}{(k+\frac{n+2}{2})} \{\sum_{m=1}^{k+1} \frac{1}{(m+\frac{n-1}{2})} t_m(n-2) + \sum_{m=1}^{k+1} \frac{1}{(k-m+\frac{3}{2})} t_m(n-2)\}
\end{aligned} \qquad (13)$$

Now use **( 6 )**, **( 8 )** and **( 10 )** to simplify, by adding and subtracting missing terms:

$$\begin{aligned}
L_{k+1}(n) &= L_k(n) + \frac{1}{(k+\frac{n+2}{2})} \{t_{k+1}(n-1) - \frac{2}{n-1} t_0(n-2) + (n-1) t_{k+1}(n-1) \\
&\quad - \frac{1}{(k+\frac{3}{2})} t_0(n-2)\} + \frac{t_0(n-1)}{(k+\frac{3}{2})} \\
&= L_k(n) + \frac{n}{(k+1+\frac{n}{2})} t_{k+1}(n-1)
\end{aligned} \qquad (14)$$

With reference to **( 11 )**, $L_{k+1}(n)=R_{k+1}(n)$. Thus the assumption is true for successive values of "$k$" for any value of "$n$", and the theorem is proven. ▲

**Corollary**

$$\sum_{m=0}^{k}\frac{1}{(k-m+1/2)(m+n/2)}t_m(n-1) = \frac{(n+1)}{k+\frac{(n+1)}{2}}\sum_{m=0}^{k}\frac{1}{(m+n/2)}t_m(n-1) \qquad (15)$$

**Proof:** Apply partial fraction decomposition to the left-hand side of **( 15 )**, using **( 9 ))**. ▲

**Theorem 3.**

The following result is the main point of this work. Consider the function $T(n,x) \equiv (\arctan(x)/x)^n$. A very simple and elegant expansion, analogous to the classic expansion of $[\log(1+x)]^n$ in terms of Stirling numbers [**1**, Section 24.1.3] is

$$(\arctan(x)/x)^n = \frac{\Gamma(n+1)}{2^n}\sum_{m=0}^{\infty}\frac{(-x^2)^m}{m+\frac{n}{2}}t_m(n-1) \qquad |x|<1 \qquad (16)$$

**Proof:** Using well-known identities, express T(n,x) as a product of logarithms, thence as a product of infinite sums, and shift the summation indices according to a well-known procedure:

$$\begin{aligned}(\arctan(x)/x)^n &= \frac{1}{2^n}\prod_{i=1}^{n}\left(\sum_{m_i=0}^{\infty}\frac{(-x^2)^{m_i}}{(m_i+1/2)}\right) \\ &= \frac{1}{2^n}\sum_{m_n=0}^{\infty}(-x^2)^{m_n}\sum_{m_{n-1}=0}^{m_n}\frac{1}{(m_n-m_{n-1}+1/2)}\cdots \\ &\cdots \sum_{m_2=0}^{m_3}\frac{1}{(m_3-m_2+1/2)}\sum_{m_1=0}^{m_2}\frac{1}{(m_2-m_1+1/2)(m_1+1/2)}\end{aligned} \qquad (17)$$

This is a known result [5, Eq. 88.2.4]. Apply **( 15 )** with (n=1) to the rightmost sum in **(17 )**, then to the next sum over $m_2$ using n=2, and propagate leftwards (n-1) times in all, each time using **( 15 )** with a unit increase in n. The result is **( 16 )**. ▲

**Corollary**: **( 16 )** can be rewritten as:

$$\left(\sum_{m=0}^{\infty} \frac{x^m}{(m+\frac{1}{2})}\right)^n = n! \sum_{0 \leq m_1 \leq m_2 \leq \cdots \leq m_n} \frac{x^{m_n}}{(m_1+\frac{1}{2})(m_2+1)\cdots(m_n+\frac{n}{2})} \qquad (18)$$

**Some properties of the coefficients** $t_m(n-1)$

**( 16 )** is equivalent to a Taylor's series expansion about x=0, (convergent for $|x| \leq 1$). Justified by this observation, equate each of the coefficients in **( 16 )** with the $(2m)^{th}$ derivative of T(n,x):

$$\lim_{x=0} \frac{d^{2m}}{dx^{2m}} T(n,x) = (-1)^m \frac{\Gamma(n+1)\Gamma(2m+1)}{2^n (m+\frac{n}{2})} t_m(n-1) \qquad m=1,2,\cdots \qquad (19)$$

The derivative calculation is straightforward, if tedious; this has been done using the Maple 6 computer code, and the first few results are given in Table 1.

Table 1. First few terms in the series expansion of T(n,x).

| m | $\lim_{x=0} \dfrac{d^{2m}}{dx^{2m}} T(n,x)$ |
|---|---|
| 1 | $-\dfrac{2}{3} n$ |
| 2 | $\dfrac{4}{15} n(5n+13)$ |
| 3 | $-\dfrac{8}{63} n(35n^2 + 273n + 502)$ |
| 4 | $\dfrac{16}{135} n(175n^3 + 2730n^2 + 13589n + 21306)$ |
| 5 | $-\dfrac{32}{99} n(385n^4 + 10010n^3 + 94259n^2 + 377938n + 538008)$ |

**( 19 )** and Table 1 generate a family of identities for a partial sum of partial sums of rational numbers. For example, the case m=3 is given below, valid for **any** n>1.

$$t_3(n-1) = \sum_{m_{n-1}=0}^{3} \frac{1}{(m_{n-1} + \frac{n-1}{2})} \sum_{m_{n-2}=0}^{m_{n-1}} \frac{1}{(m_{n-2} + \frac{n-2}{2})} \cdots \sum_{m_2=0}^{m_3} \frac{1}{(m_2 + 1)} \sum_{m_1=0}^{m_2} \frac{1}{(m_1 + \frac{1}{2})}$$
$$= \frac{2^{n+3}(3 + n/2)}{63(n-1)!\,6!}(35n^2 + 273n + 502)$$
( 20 )

Finally, using the easily derived property that

$$\frac{x^2}{n} \frac{d^2}{dx^2} T(n,x) = (n+1)T(n,x) - \frac{2(n + (n+1)x^2)}{(1+x^2)^2} T(n-1,x) + \frac{(n-1)}{(1+x^2)^2} T(n-2,x) \quad (21)$$

it is possible to find a 5-part recursion formula for the coefficient $t_m(n)$ by substituting **( 16 )** into **( 21 )** and equating coefficients of equal powers of x:

$$t_m(n) = \Big(-2(n - 2m + 4)t_{m-1}(n) + (n - 2m + 6)t_{m-2}(n) - \frac{4(n+1)}{2m+n} t_m(n-1)$$
$$+ \frac{4(n+2)}{(2m+n-2)} t_{m-1}(n-1) + \frac{4}{2m+n-1} t_m(n-2)\Big)/(2m - n - 2)$$
( 22 )

This result is valid if $(2m - n - 2) \neq 0$ in which case the numerator also vanishes, (implying a four-part recursion for this special case) and another method must be used (e.g. **( 31 )**).

In terms of digamma functions, the first two known [**5**, Eqs. 55.6.3 and 5.4.26], specific values are:

$$t_m(1) = \psi(m+\tfrac{3}{2}) - \psi(\tfrac{1}{2})$$

$$t_m(2) = \sum_{k=0}^{m} \frac{\psi(k+\tfrac{3}{2}) - \psi(\tfrac{1}{2})}{(k+1)} \qquad (23)$$

From **(16)** a number of interesting sums can be found, corresponding to values of "$x$" for which arctan has a known value ($x = 2-\sqrt{3}$, $1-\tfrac{2}{\sqrt{5}}$, etc.). For example, consider the barely convergent sum **(16)** at $x=1$

$$\frac{\pi^n}{2^n \Gamma(n+1)} = \sum_{m=0}^{\infty} \frac{(-1)^m}{m + \tfrac{n}{2}} t_m(n-1) \qquad (24)$$

which reduces to a well-known classical result when n=1.

**Lemma**

$$i^m \frac{c}{2i} \sum_{l=0}^{m} \frac{\Gamma(m + \tfrac{c}{2i} - l)}{\Gamma(1+m-l)\Gamma(1+\tfrac{c}{2i}-l)\Gamma(l+1)} = \frac{2}{m} \sum_{l=0}^{\left[\tfrac{m}{2}\right]} (\tfrac{c}{2})^{m-2l} (-1)^l t_l(m-2l-1) \qquad (25)$$

where $\left[\tfrac{m}{2}\right]$ is the greatest integer less than or equal to m/2, m>0 and "$c$" is continuous.

**Proof:** Consider the sum

$$\exp(c \, \arctan(x)) = \sum_{n=0}^{\infty} \frac{c^n (\arctan(x))^n}{\Gamma(n+1)} = \sum_{m=0}^{\infty} (-x^2)^m \sum_{n=0}^{\infty} \left(\frac{cx}{2}\right)^n t_m(n-1)/(m+\tfrac{n}{2}) \qquad (26)$$

But

$$\exp(c \, \arctan(x)) = \left(\frac{1+ix}{1-ix}\right)^{\tfrac{c}{2i}} = i^m \sum_{m=0}^{\infty} g_m(\tfrac{c}{2i}) x^m \qquad (27)$$

using the well-known identification $\arctan(x) = \log(\frac{1+ix}{1-ix})/2i$. The coefficients $g_m(\frac{c}{2i})$ in **( 27 )** are identified [**4**, Eq. 19.6(23)] as

$$g_m(y) = \frac{(y)_m}{m!} {}_2F_1(-m,-y;1-m-y;-1) \qquad (28)$$

where $(y)_m$ is Pochhammer's symbol. Write the terminating hypergeometric function as an explicit sum, let y=c/(2i) and simplify. ▲

**Theorem 4.**

$$\sum_{l=0}^{\left[\frac{m}{2}\right]} t_l(m-2l-1) = m \qquad (29)$$

**Proof:** This sum rule is proven by setting c=2i in **( 25 )**. The terminating sum has only two terms, and the proof is immediate. ▲

**Theorem 5.**
Since the left-hand side of **( 25 )** reduces to a polynomial of degree "m-1" in "c", emphasized by writing the "c" dependence using the Pochhammer symbol $(1+\frac{c}{2i}-l)_{m-1}$, let

$$i^m \frac{c}{2i} \sum_{l=0}^{m} \frac{(1+\frac{c}{2i}-l)_{m-1}}{\Gamma(1+m-l)\Gamma(l+1)} = \sum_{l=0}^{\left[\frac{m}{2}\right]} (\frac{c}{2i})^{m-2l} p_l(m) \qquad (30)$$

Then

$$t_l(n-1) = \frac{m}{2}(-1)^l p_n(m) \quad \text{where} \quad m = n+2l \qquad (31)$$

**Proof:** Since "c" is arbitrary, equate coefficients of equal powers of "c" in **( 25 )** and **( 30 )**. ▲

From **( 25 )** it is now possible to obtain a closed expression for the coefficients $p_n(m)$ and hence reduce the n-fold sum, symbolized by $t_l(n-1)$ to a simpler form.

**Theorem 6.**

$$p_n(m) = \frac{1}{\Gamma(n)} \sum_{j=0}^{m-n} \frac{\Gamma(j+n)}{\Gamma(j+1)} L_m^j S_{m-1}^{j+n-1} \qquad ( 32 )$$

where $S_{m-1}^{j+n-1}$ are Stirling numbers of the first kind [**1**, Section 24.1.3] and $L_m^j$, given in terms of generalized Bernoulli polynomials [**5**, Eq.6.7.26], are coefficients defined by

$$L_m^j = \sum_{l=0}^{m} \frac{(l-1)^j}{\Gamma(1+l)\Gamma(1+m-l)} \qquad ( 33 )$$

In **( 33 )**, $(l-1)^j = 1$ when $l = 1$ and $j = 0$.

**Proof:** On the left-hand side of **( 30 )**, reverse the series, and expand as follows:

$$\frac{\Gamma(C+l)}{\Gamma(C+l+1-m)} = \sum_{j=0}^{m-1} S_{m-1}^j (C+l-1)^j = \sum_{j=0}^{m-1} S_{m-1}^j \sum_{k=0}^{j} C^k (l-1)^{j-k} \binom{j}{k} \qquad ( 34 )$$

using well-known expansions [**4**, Eq. 19.7(61), **1**, Eq. (24.1.3B)], and setting $C = c/(2i)$. Invert and re-order the double sum. ▲

**References**

**1** Abramowitz,M. and Stegun,I.,"Handbook of Mathematical Functions", National Bureau of Standards (1964).

**2.** Bromwich, T.J., "An Introduction to the Theory of Infinite Series", Macmillan and Co., (1908, revised (1926)).

**3.** Case, K.M., de Hoffman, F and Placzek, G. "Introduction to the Theory of Neutron Diffusion", LosAlamos report, (1953).

**4.** Erdelyi,A.,Magnus,W., Oberhettinger, F., Tricomi, F.G., "Higher Transcendental Functions", Vol. III, McGraw-Hill (1955).

**5.** Hansen, E.R, "A Table of Series and Products", Prentice-Hall (1975).

**6.** Milgram, M.S., "The Generalized Integro-Exponential Function", Math. Comp. 44,170, pp 443-458 (1985).

**7.** Milgram, M.S. "On Some Sums of Digamma and Polygamma Functions", submitted for publication, available as **http://www.arXiv.org: math:CA/0406338, (2004a).**

**8.** Milgram, M.S. "Identification and Properties of the Fundamental Expansion Functions for Neutron Transport in a Homogeneous Scattering Medium", (in preparation), (2004b).

**9** . Prudnikov, A.P., Brychkov, Yu.A., Marichev, O.I., "Integrals and Series", Vol.3, "More Special Functions", Gordon and Breach, 1990.

**10.** Schwatt, I.J., "An Introduction to the Operations with Series", University of Pennsylvania Press, (1924).

**11.** Wolfram, (2004): http://functions.wolfram.com/